\documentclass[12pt]{elsarticle}
\usepackage{setspace}
\usepackage{amsmath,amscd, amsthm}
\usepackage{amsfonts}
\usepackage[english]{babel}
\usepackage{amssymb}
\usepackage[top=1in, bottom=1in, left=1in, right=1in]{geometry}
\usepackage{enumerate}
\usepackage{mathrsfs}

\usepackage{xypic}
\input xy

\newcommand{\R}{\mathbb{R}} 
 
\newcommand{\C}{\mathbb{C}}

\newcommand{\tab}{\hspace*{15pt}}

\newcommand{\xqedhere}[2]{%
  \rlap{\hbox to#1{\hfil\llap{\ensuremath{#2}}}}}

\usepackage{tikz-cd}
\usepackage{enumerate}
\usepackage{mathrsfs}
\usepackage[utf8]{inputenc}
\usepackage[T1]{fontenc}
\usepackage{xypic}
\input xy

\newtheorem{thm}{Theorem}

\newtheorem{lem}[thm]{Lemma}

\xyoption{all}

\title{A spectral interpretation of zeros of certain functions}
\author{Kim Klinger-Logan}
\date{01.27.2020} 

\begin{document}

\maketitle

\begin{quote}{\sc Abstract:} We prove that all the zeros of certain meromorphic functions are on the critical line $\text{Re}(s)=1/2$, and are simple (except possibly when $s=1/2$). We prove this by relating the zeros to the discrete spectrum of an unbounded self-adjoint operator. Specifically, we show for $h(s)$ a meromorphic function with no zeros in $\text{Re}(s)>1/2$ and no poles in $\text{Re}(s)<1/2$, real-valued on $\R$, $\frac{h(1-s)}{h(s)}\ll |s|^{1-\epsilon}$ in $\text{Re}(s)>1/2$ and $\frac{h(1-s)}{h(s)}\notin L^2(1/2+i\R)$, the only zeros of $h(s)\pm h(1-s)$ are on the critical line.  One instance of such a function $h$ is $h(s)=\xi(2s)$, the completed zeta-function. 
We use spectral theory suggested by results of Lax-Phillips and Colin de Verdi\`{e}re.  This simplifies ideas of W. M\"{u}ller, J. Lagarias, M. Suzuki, H. Ki,  O. Vel\'{a}squez Casta\~{n}\'{o}n, D. Hejhal, L. de Branges  and P.R. Taylor.
\\ \end{quote}

One method of showing that the zeros of a function lie on the critical line $1/2+i\R$ is to identify those zeros with spectral parameters of a self-adjoint operator. 
The simplest unbounded operator that may be used in this context is a multiplication operator on a Hilbert space.  In what follows, we provide a set of conditions for which the zeros of a function may appear as spectral parameters for a multiplication operator on a Hilbert space.

\tab In order to accurately state our result, we introduce the following Hilbert spaces of functions.  
Let $\eta=\pm 1$ and 
$$H_\eta^0:=\left\{f\in L^2(1/2+i\mathbb{R})~\Big|~f(1-s)=\eta\cdot \frac{f(s)}{c_s} \text{ for Re}(s)=\frac{1}{2} \right\}$$ with $L^2$-norm and put $D:=H_\eta^0\cap {C_c^\infty(1/2+i\mathbb{R})}$.\footnote{ Note that for what follows we could use the full space $L^2(1/2+i\mathbb{R})$ without restricting to those which satisfy the functional equation. This larger space will still be stable under the multiplication operator we define and the proof will follow in a similar way. In this paper, we restrict our attention to the smaller space $H_\eta^0$ as in subsequent work and applications of this technique there are reasons for such a restriction.}
 For each $ k\in\mathbb{R}$, let $H_\eta^k$ be the Sobolev-like space
$$H_\eta^k:=\Big\{f\text{ measurable on }1/2+i\mathbb{R}~\Big|~\langle \left (s(1-s)\right)^{k}\cdot f, f\rangle _{L^2(1/2+i\mathbb{R})}<\infty\phantom{weeeeeeeee}$$
 \vspace{-.7cm}
$$\phantom{weeweeeeweeeeeeeeeeeeeeweeeeeeeeeeeeeeeeeeeeeeee}\text{ and } f(1-s)=\eta  \cdot \frac{f(s)}{c_s} \text{ for Re}(s)=\frac{1}{2}  \Big\}$$ with norm-squared $||f||^2_{H_\eta ^k}=\langle \left(s(1-s)\right)^{k}f, f \rangle_{L^2}.$
Lemma \ref{duals} shows that  $H ^1_\eta$ and $H ^{-1}_\eta$ are mutual Hilbert-space duals and functions in $H ^{-1}_\eta$ give functionals on $H ^1_\eta$. 

\tab In this paper we will prove the following:

{\thm \label{nozeros}Let $h(s)$ be a meromorphic function with no zeros in $\text{Re}(s)>1/2$, no poles in $\text{Re}(s)<1/2$, $\mathbb{R}$-valued on $\R$, and $c_s:=\frac{h(1-s)}{h(s)}\ll |s|^{1-\epsilon}$ uniformly in $\text{Re}(s)>1/2$ for some $\epsilon >0$. For $\eta=\pm 1$, assume $1+ \eta\cdot  c_s$ is in $H_\eta^{-1}\setminus H^0_\eta$.
Then the only zeros of $1+ \eta \cdot c_s$ are on the critical line. 
 Furthermore, all of the zeros on the critical line are simple (with the possible exception of $s_0=1/2$).}\\

\tab Some examples of functions $h$ that satisfy the hypotheses of Theorem \ref{nozeros} are $h(s) =\xi(2s)$ where $\xi$ is the completed Riemann-zeta function, $h(s)=\xi_k(2s)$ the completed zeta-function of a number field $k$ and many self-dual automorphic $L$-functions, and $h(s)=L(2s,\chi)$ a Dirichlet $L$-function. Theorem \ref{nozeros} shows that $\xi(2s)\pm \xi(2-2s)$ has all of its zeros on the critical line $1/2+i\R$. On of the first examples of a result of this kind is due to P.R. Taylor who proved a similar conclusion for $h(s)=\xi(s-1/2)$  \cite{Taylor}.  
Results of a similar form have been established for various functions $h$ by Lagarias-Suzuki \cite{Lagarias, Lagarias2006}, Ki \cite{Ki}, McPhedran-Poulton \cite{McPhadren}, Vel\'{a}squez Casta\~{n}\'{o}n \cite{Velasquez},  Hejhal \cite{Hejhal1990}, and Taylor \cite{Taylor}. 
We use more general methods, vaguely reminiscent of de Branges \cite{deBranges}, and expanded-upon by Kaltenb\"{a}ck and Woracek \cite{Kaltenback2004}.  The Lax-Phillips 1976 \cite{LP1976} automorphic example arguably suggests a similar result for constant terms of Eisenstein series on reductive groups, and M\"{u}ller's work \cite{Muller} is an extension of Lax-Phillips. The Lax-Phillips-ColinDeVerdi\`{e}re-Hejhal-M\"{u}ller \cite{LP1976, Hejhal1990, Muller} base example was $h(s)=\xi(2s)$.  Lagarias-Suzuki \cite{Lagarias2006}, Hejhal  \cite{Hejhal1990} and M\"{u}ller \cite{Muller} also treated the case $h(s)=\xi(2s)y^s$ to show that all the zeros of the constant term of the Eisenstein series are on the critical line. 
 
 \tab Small modifications must be made to our result to subsume the above examples. For instance, the construction we present in the proof of Theorem \ref{nozeros} may be extended to include the case where $h$ has a finite number of poles on the real line.

\tab It is worth noting that, though Taylor's paper suggests that some have hoped otherwise, this approach cannot prove that the nontrivial zeros of $\zeta (s)$ are on the critical line.  In fact, any means applying the argument we outline to genuine zeta functions would similarly apply to Epstein zeta functions which are known to have many off-line zeros (see Potter-Titchmarsh \cite{PT1935}, Stark  \cite{Stark1967}, Voronin \cite{Voronin1976}, et al).

\tab  In Section 1, we prove that $H^k_\eta$ and $H^{-k}_\eta$ defined above are Hilbert space duals. In Section 2, we explain that the proof of Theorem {\ref{nozeros}} is given by identifying the zeros of the function $1\pm c_s$ with the spectral parameters for eigenvalues $\lambda_s=s(1-s)$ of the Friedrichs extension of an unbounded operator.  
As explained in the Appendix (Section 4), the Friedrichs extension is a self-adjoint extension and so this construction will show that all of the eigenvalues are real. In Section 3, we will prove the main part of Theorem 1 (that all the zeros of $1+\eta\cdot c_s$ lie on the critical line). Finally, in Section 4 we will show that these zeros are simple.


\section{A note on Hilbert spaces}

Let $H_\eta^0$ and $H_\eta^k$ be as defined in the introduction.  The {\it complex bilinear} paring putting $H_\eta^k$ and $H_\eta^{-k}$ in (complex-linear duality) is 
$$\langle f, F \rangle_{H_\eta^k\times H_\eta^{-k}} = \int_{(1/2)}f(s)\cdot F(1-s)\, ds 
=\int_{(1/2)}f(s)\cdot\frac{\eta}{c_s}\cdot  F(s)\, ds$$ since $F\in H^{-k}$ requires that $F(1-s)= \frac{\eta}{c_s}\cdot F(s)$. Let $F^{\vee}(s) = F(1-s)$
where we use the notation $\displaystyle{\int_{(1/2)}g(s)~ds:=\int_{1/2-i\infty}^{1/2+i\infty}g(s)~ds }$. Then the {\it hermitian} pairing (conjugate-linear in the second argument) on $H^{k}\times H^{-k}$ is 
$$\langle f, F\rangle_{herm} = \langle f, \overline{F^{\vee}}\rangle_{H_\eta^{k}\times H_\eta^{-k}} = \int_{(1/2)}f(s)\cdot \overline{F(s)}\, ds$$

{\lem\label{duals}For $k\in\R$, $H_\eta^k$ and $H_\eta^{-k}$ are mutual $\C$-linear Hilbert-space duals by the complex-bilinear pairing $\displaystyle \langle g,f\rangle_{H_\eta^{k}\times H_\eta^{-k}}= \int_{(1/2)}g(s)\cdot { f(1-s)}\,ds$.}

\begin{proof}  Let $f\in H_\eta^{-k}$ and define $\varphi_f: H_\eta^k\to \mathbb{C}$ by $\varphi_f(g)=\langle g,f\rangle_{H_\eta^{k}\times H_\eta^{-k}}$ for each $g\in H_\eta^k$. Let $\lambda_s=s(1-s)$ and $f\in H_\eta^{-k}$ and $g\in H_\eta^k$
\begin{align*}
| \varphi_f(g)|&=\left|\langle g,f\rangle_{H_\eta^{k}\times H_\eta^{-k}}\right| 
 = \left|\int_{(1/2)}g(s)\cdot  { f(1-s)}\,ds\right|
 \leq \int_{(1/2)}\left| g(s)\cdot { f(1-s)}\right| \,ds\\
&= \int_{(1/2)}\left|  (1+ |\lambda_s|)^{k/2}\cdot g(s)\cdot (1+ |\lambda_s|)^{-k/2}\cdot  { f(1-s)}\right| \,ds\\
&\leq \left(\int_{(1/2)}\left|g (s)\right|^2 (1+ |\lambda_s|)^{k} \,ds\right)^{1/2}\cdot \left(\int_{(1/2)}\left| f(1-s)\right|^2 (1+ |\lambda_s|)^{-k} \,ds\right)^{1/2}
\end{align*}
  \hfill by the  Cauchy-Bunyakovsky-Schwarz Inequality
  $$= |g|_{H_\eta^{k}}\cdot|f|_{H_\eta^{-k}}<\infty$$
  Since $\varphi_f$ is bounded, it is continuous.
  
 \tab  Conversely by Riesz-Fr\'{e}chet, every continuous functional on $H_\eta^0$ is given by integration against an element of $H_\eta^0$. Given $F\in H_\eta^k$, we have $F(s)\cdot (1+ |\lambda_s|)^{k/2}\in H_\eta^0$ and $F\mapsto F(s)\cdot (1+ |\lambda_s|)^{k/2}$ is an isomorphism $H_\eta^k\to H_\eta^0$. Thus every continuous linear functional $\rho$ of $H_\eta^k$ factors through $H_\eta^0$:
 $$\rho(F)= \int_{(1/2)} F(s)\lambda_s^{k/2}\cdot {g(s)}\,ds$$
 for $g\in H_\eta^0$. Similarly, $g(s) ={G}(1-s) (1+ |\lambda_s|)^{-k/2}$ for unique $G\in H_\eta^{-k}$ and 
 $$\rho(F)= \int_{(1/2)} F(s)\lambda_s^{k/2}\cdot {G(1-s)} (1+ |\lambda_s|)^{-k/2}\,ds =\int_{(1/2)} F(s)\cdot {G(1-s)}\,ds$$
 This gives the duality pairing $H_\eta^k\times H_\eta^{-k}\to \C$.
 
 \end{proof}

Thus, functions in $H_\eta^{-1}$ give functionals on $H_\eta^1$.

\section{Overview}\label{exp}

\tab The idea behind the proof of Theorem \ref{nozeros} is to identify the zeros of a function $1\pm c_s$ with spectral parameters $s$ for eigenvalues $\lambda_s=s(1-s)$ of the Friedrichs extension $\widetilde T$ of $T$.
The Friedrichs extension (discussed in generality in the Appendix) is a self-adjoint extension of a densely-defined, semi-bounded, symmetric operator $T$ with domain $D$.  
The symmetry of $\widetilde T$ gives $\lambda_s\in \mathbb{R}$, from which $\text{Re}(s)=1/2$ for $T\geq 1/4$ since $\widetilde T$ has the same lower bound.
 It may be helpful to note that the Friedrichs extension of the multiplication operator is constructed in essentially the same way as the self-adjoint extension of the multiplication operator by the independent variable in the de Branges space. To construct this space we could take $c_s$ to be the Hermite-Biehler function; however, in our case $c_s$ need not be entire.

\tab Before the proof of Theorem \ref{nozeros}, we review the characterization of the Friedrichs extension and explain how the argument presented in Section \ref{proof} serves as a proof. 
Suppose $H^0$ is a Hilbert space with norm $\langle\cdot, \cdot \rangle$.  Begin with a densely-defined, semi-bounded, symmetric operator $T$ with domain $D = H^0\cap C_c^\infty$.  We can extend $T$ to a Hilbert space $H^1\supset D$ with norm $\langle u, v \rangle_1 = \langle Tu, v \rangle$ by letting { $T ^{\#}: H ^1\to H ^{-1} $ }be the continuous $\C$-linear operator defined by $$\left(\left( T ^{\#}\right)f\right)(g)=\langle f,\overline{g}\rangle_{H ^1}$$ for all $f,g\in H ^1$. 

\tab For $\theta \in H ^{-1}\setminus H^0$, let $\widetilde\theta : H ^1\to\mathbb{C}$ be the functional on $H ^1$ given by the function $\theta \in H ^{-1}$. Define $T_{\theta }:=T ^\#|_{\text{ker}\,\widetilde\theta \,\cap D}$ and let $\widetilde T_{\theta }$ be the Friedrichs extension of $T_{\theta }$.   A useful characterization of the Friedrichs extension is that for  $u\in H ^1$ and $f\in H ^0$ 
$$\widetilde T_{\theta }\, u = f \phantom{weeee}\Longleftrightarrow\phantom{weeee} T ^{\#}u = f+\alpha\cdot \theta \text{ for some }\alpha\in\C \phantom{wee} \& \phantom{wee}  \widetilde\theta  u = 0$$
as in the Appendix.  Theorem \ref{Fried} establishes that, for $u\in H ^1$,
\begin{equation*}\label{iff}(\widetilde T_{\theta }-\lambda_w)u=0 \phantom{wee} \Longleftrightarrow\phantom{we}  (T ^\#-\lambda_{w})u=\alpha\cdot \theta \phantom{wee}\text{ for some } \alpha\in\mathbb{C}\phantom{we}\& \phantom{we}\widetilde\theta ( u)=0\end{equation*}
where $\lambda_w= w(1-w).$ We may assume that $\alpha\neq 0$ since otherwise $u$ would be an eigenfunction for $T$ and there are none. Furthermore, the value of $\alpha$ does not affect the location of the zeros, so we may assume $\alpha=1$.

\tab To prove Theorems \ref{nozeros}:\vspace{-.3cm}

\begin{quote}

\begin{enumerate}[(I)]
\item Identify an appropriate $\theta $ such that the zeros $w$ of $1+  c_s$ 
$$(T^\#-\lambda_{w})u= \theta $$
with $ \lambda_w=w(1-w).$
\item Solve $(T^\#-\lambda_{w})u_w=\theta $ for $u_w$.
\item Show that $\widetilde\theta (u_w)=0$.

\end{enumerate}\end{quote}\vspace{-.2cm}

This will ensure that the zeros $w$ of $1+  c_s$ appear as parameter-values for $\lambda_w$ which satisfy  $(\widetilde T_{\theta }-\lambda_w)u=0$ where $\widetilde T_{\theta }$ is the (self-adjoint) Friedrichs extension of $T_{\theta }=T|_{\text{ker}\widetilde{\theta }}$.  Since $\widetilde T_{\theta }\geq 1/4$,  it follows that $\lambda_w=w(1-w)$ must be real and $\lambda_w\geq  1/4$ and so $\text{Re}(w)=1/2$.

\section{Proof of Theorem \ref{nozeros}}\label{proof} 

First recall the assumptions of Theorem \ref{nozeros}. Let $h(s)$ be a meromorphic function with no zeros in $\text{Re}(s)>1/2$, no poles in $\text{Re}(s)<1/2$, $\mathbb{R}$-valued on $\R$, and $c_s:=\frac{h(1-s)}{h(s)}\ll |s|^{1-\epsilon}$ uniformly in $\text{Re}(s)>1/2$, for some $\epsilon >0$. For $\eta=\pm 1$, assume $1+ \eta\cdot  c_s$ restricted to the critical line $1/2+i\mathbb{R}$ is not square-integrable.

\tab For the proof of Theorem \ref{nozeros}, we first claim that it is sufficient to show that $1+\eta\cdot c_s$ has no zeros in $Re(s)>1/2$. 

{\lem{If $1+\eta\cdot c_s$ has no zeros in $Re(s)>1/2$ then it also has no zeros in $Re(s)<1/2$.}}
\begin{proof}Since $c_sc_{1-s}=\frac{h(1-s)}{h(s)}\cdot \frac{h(s)}{h(1-s)}= 1$, we have $\eta\cdot c_s(1+\eta\cdot c_{1-s})=1+\eta\cdot c_s$.  The numerator $h(1-s)$ of $c_s$ has no zeros in $\text{Re}(s)<1/2$ since $h(s)$ has no zeros in $\text{Re}(s)>1/2$. If the denominator $h(s)$ had a pole at $s_o$ in $\text{Re}(s)<1/2$ (causing a zero of $c_s$), there would also be a pole of $c_{1-s}=\frac{h(s)}{h(1-s)}$ (and of $1+ \eta\cdot c_{1-s}$) at the same point which will cancel. Thus  we have that $\eta\cdot c_s(1+\eta\cdot c_{1-s})=1+\eta\cdot c_s$ is nonzero in $\text{Re}(s)<1/2$, proving the claim. 
\end{proof}

\tab Define $\theta_\eta$ is $\theta_\eta(s):=1+\eta\cdot c_s$.  Let $\widetilde\theta_\eta: H_\eta^1\to \mathbb{C}$ be the functional on $H_\eta^1$ corresponding to the function $\theta_\eta\in H_\eta^{-1}$: for $f\in H_\eta^1$,

 $$\widetilde\theta_\eta(f)=\int_{(1/2)}{\theta_\eta(s) }\cdot{ f(1-s)}\,ds.$$ Lemma \ref{duals} shows that $\widetilde\theta_\eta\in (H_\eta^{1})^*$ the complex-bilinear dual of $H_\eta^{1}$. Let $T$ be multiplication operator by $s(1-s)$ on $1/2+i\R$ with domain $D=H_\eta^0\cap {C_c^\infty(1/2+i\mathbb{R})}$.
We can extend $T$ to $T_\eta^\#:H_\eta^1\to H_\eta^{-1}$ as above by letting { $T_\eta^{\#}: H_\eta^1\to H_\eta^{-1} $ } be the continuous $\C$-linear operator defined by $$\left(\left(  T_\eta^{\#}\right)f\right)(g)=\langle f,\overline{g}\rangle_{H_\eta^1}$$ for all $f,g\in H_\eta^1$. Observe that, in this case, $T_\eta^{\#}$ is still a kind of multiplication operator but mapping $H^1_\eta$ to $H^{-1}_\eta$ rather than $D$ to $H^0_\eta$.

\tab In Lemma \ref{sol1} we find a solution $u$ to $(T_\eta^\#-\lambda_{w})u= \theta_\eta$ for $\text{Re}(w)>1/2$.  

{\lem\label{sol1}  For $\text{Re}(w)>1/2$, there is a solution $u_w\in H_\eta^1$ to $(T_\eta^\#-\lambda_{w})u_w= \theta_\eta$.}

\begin{proof} Solving the equation by division, $$u_w(s)=\frac{ 1+\eta\cdot c_s}{\lambda_s-\lambda_w}\phantom{weeee} \text{on Re}(s)=1/2\text{.}$$ 
We check that this is a solution by seeing that, on $\text{Re}(s)=1/2$,
$$(T_\eta^\#-\lambda_{w})u_w = (T_\eta^\#-\lambda_{w})\frac{ 1+ c_s}{\lambda_s-\lambda_w}
=\frac{ (T_\eta^\#-\lambda_{w})(1+ \eta\cdot c_s)}{\lambda_s-\lambda_w} = \frac{ (\lambda_s-\lambda_{w})(1+ \eta\cdot c_s)}{\lambda_s-\lambda_w} = \theta_\eta.$$

\tab To see that $u_w\in H_\eta^1$, note that it satisfies the functional equation on $\text{Re}(s)=1/2$: 
$$\eta\cdot c_su_w(1-s) =\eta\cdot c_s\cdot\frac{1+\eta\cdot c_{1-s}}{\lambda_{1-s}-\lambda_w} = \frac{\eta\cdot c_s+1}{\lambda_{s}-\lambda_w} = u_w(s)$$
\hfill since $c_sc_{1-s}=1$ and $\lambda_s=\lambda_{1-s}$.  

Furthermore, 

$$\langle u_w, u_w \rangle_{H_\eta^1} = 
\int_{(1/2)}\lambda_s\cdot u_w(s)\cdot \overline{u_w(s)}\,ds
=\int_{(1/2)}\lambda_s\cdot\frac{1+\eta\cdot c_s}{\lambda_s-\lambda_w}\cdot \overline{\frac{1+\eta\cdot c_s}{\lambda_s-\lambda_w}}\,ds$$
$$=\int_{(1/2)}(1+  \eta\cdot c_s)(1+\eta\cdot c_{1-s})\frac{\lambda_s}{(\lambda_s-\lambda_w)^2}\,ds
=\int_{(1/2)}(1+ \eta\cdot c_s)(1+\eta\cdot c_{1-s})\frac{\lambda_s^2}{(\lambda_s-\lambda_w)^2\cdot \lambda_s}\,ds$$

Since $\theta_\eta\in H_\eta^{-1}$, we have that $\theta_\eta(\overline{s})\frac{\lambda_s^2}{(\lambda_s-\lambda_w)^2}=(1+ \eta\cdot c_{1-s})\frac{\lambda_s^2}{(\lambda_s-\lambda_w)^2}\in H_\eta^{-1}$. 
Thus the expansion above is equal to 
$$ \left\langle \theta_\eta(s), \theta_\eta(\overline{s})\frac{\lambda_s^2}{(\lambda_s-\lambda_w)^2}\right \rangle_{H_\eta^{-1}}$$ which is finite since $\theta_\eta\in H_\eta^{-1}$ and we see that $u_w\in H_\eta^1$.

\end{proof}

{\lem\label{BC1} For $\text{Re}(w)>1/2$, if $\theta_\eta(w)=  0$ then $\widetilde\theta_\eta(u_w)= 0$.}

\begin{proof}

Assume $\text{Re}(w)>1/2$. We have
$$\widetilde\theta_\eta(u_w)=\int_{(1/2)}{(1 + \eta\cdot  c_s)} u_w(1-s)\,ds
\phantom{www}\text{for Re}(w)>1/2.$$

\tab We compute the integral $\displaystyle{\int_{(1/2)}{(1 +\eta\cdot c_s)} u_w(1-s)\,ds}$ as follows:
$$\int_{(1/2)}{(1 + c_s)} u_w(1-s)\,ds=
\int_{(1/2)} \frac{{(1+\eta\cdot c_s)}\cdot(1+\eta\cdot  c_{1-s})}{\lambda_s-\lambda_w}\,ds $$
    $$=\int_{(1/2)} \frac{1+ \eta\cdot c_s+\eta\cdot c_{1-s}+c_s{c_{1-s}}}{\lambda_s-\lambda_w}\,ds$$
\hfill since $\displaystyle{u_w(s)=\frac{1+ \eta\cdot c_s}{\lambda_s-\lambda_w}}$ for $s\in 1/2+i\R$
$$=\int_{(1/2)} \frac{2}{\lambda_s-\lambda_w}\,ds +\eta \int_{(1/2)} \frac{c_s}{\lambda_s-\lambda_w}\,ds+\eta \int_{(1/2)} \frac{c_{1-s}}{\lambda_s-\lambda_w}\,ds$$
$$=2\int_{(1/2)} \frac{1}{\lambda_s-\lambda_w}\,ds + 2\eta \int_{(1/2)} \frac{c_s}{\lambda_s-\lambda_w}\,ds$$ \hfill by a change of variables $s\mapsto 1-s$.\\

For each of these integrals (with corresponding function $g(s)$), we have
$$\int_{(1/2)}\frac{g(s)}{\lambda_s-\lambda_w}\,ds=\lim_{T\to \infty}\int_{1/2-iT} ^{1/2+iT}\frac{g(s)}{\lambda_s-\lambda_w}\,ds.$$
For any $\sigma>\sigma_n$, let $S_T$ be a clockwise oriented closed semi-circle of radius $T$ with endpoints at $1/2-iT$ and $1/2+iT$ and let $A_T$ be the outer arc of $S$ with length $\pi T$. $$\int_{1/2-iT} ^{1/2+iT}\frac{g(s)}{\lambda_s-\lambda_w}\,ds=\int_{S_T} \frac{g(s)}{\lambda_s-\lambda_w}\,ds-\int_{A_T}  \frac{g(s)}{\lambda_s-\lambda_w}\,ds$$ \\
As we send $T\to \infty$ each of the integrals $\int_{A_T} \frac{g(s)}{\lambda_s-\lambda_w}\,ds\to 0$:
$$\left|\int_{A_T} \frac{1}{\lambda_s-\lambda_w}\,ds\right|\leq\int_{A_T}\frac{1}{|\lambda_s-\lambda_w|}\,ds= \int_{A_T} \frac{1}{|s^2-s-\lambda_w|}\,ds\leq\frac{1}{|T^2-T-\lambda_w|}\cdot\pi T\leq\pi T^{-1}$$
Since we assume $c_s:=\frac{h(1-s)}{h(s)}\ll |s|^{1-\epsilon}$ uniformly in $\text{Re}(s)>1/2$ for some $\epsilon >0$, we also have
$$\left|\int_{A_T} \frac{c_s}{\lambda_s-\lambda_w}\,ds\right|\leq\int_{A_T} \frac{|c_s|}{|\lambda_s-\lambda_w|}\,ds\leq\int_{A_T} \frac{|s|^{1-\epsilon}}{|s^2-s-\lambda_w|}\,ds
\leq T^{-1-\epsilon}\cdot\pi T=\pi T^{-\epsilon}$$ for $0<\epsilon<1$, each of these integrals approaches $0$ as $T\to \infty$.

\tab Thus we can use the Residue Theorem to compute $$\int_{(1/2)} \frac{g(s)}{\lambda_s-\lambda_w}\,ds.$$ Observe that $c_s$ has no poles on $\text{Re}(s)=1/2$ since $|h(s)|=|h(1-s)|$ on $\text{Re}(s)=1/2$.  Recall that $h(s)$ has no poles in $\text{Re}(s)<1/2$, so $h(1-s)$ has no poles in $\text{Re}(s)>1/2$. Also, $h(s)$ no zeros in $\text{Re}(s)>1/2$ so $c_s$ has no poles in $\text{Re}(s)>1/2$ and $ \frac{g(s)}{\lambda_s-\lambda_w}$ has a simple pole at $s=w$. By residues, for $\text{Re}(w)>1/2$
 $$2\alpha\int_{(1/2)} \frac{1}{\lambda_s-\lambda_w}\,ds+2\alpha\cdot \eta \int_{(1/2)} \frac{c_s}{\lambda_s-\lambda_w}\,ds
  =-2\pi  i\left(2\,\text{Res}_{s=w}\frac{1}{\lambda_s-\lambda_w} +2\eta\cdot \text{Res}_{s=w}\frac{c_s}{\lambda_s-\lambda_w}\right)$$
 $$ =-2\pi i\left(\frac{2}{1-2w} + \frac{2\eta\cdot c_w}{1-2w} \right)=-4\pi i\left( \frac{1+\eta\cdot c_w}{1-2w} \right)$$
 \hfill since $h(w)$ has no zeros in $\text{Re}(s)>1/2$

 \tab We then have $$\widetilde\theta_\eta(u_w)= -4\pi i\left( \frac{1+\eta\cdot c_w}{1-2w} \right) =-4\pi i \frac{\theta_\eta(w)}{2w-1}=0$$
\end{proof}

\tab Define $T_{\theta_\eta}:=T_\eta^\#|_{\text{ker}\,\widetilde\theta_\eta\,\cap D}$ and let $\widetilde T_{\theta_\eta}$ be the Friedrichs extension of $T_{\theta_\eta}$. 
From Lemmas \ref{sol1} and \ref{BC1}, if $\text{Re}(w_o)>1/2$, then 
$(T_\eta ^\#-\lambda_{w_o})u_{w_o}= \theta_\eta$ has a solution in $H_\eta^1$ and if $\theta(w_o)=0$ then $\widetilde\theta ( u_{w_o})=0$.
Since
\begin{equation*}(\widetilde T_{\theta }-\lambda_w)u=0 \phantom{wee} \Longleftrightarrow\phantom{we}  (T_\eta ^\#-\lambda_{w_o})u_{w_o}= \theta_\eta \phantom{wee}\& \phantom{we}\widetilde\theta ( u_{w_o})=0\end{equation*}
 for $\text{Re}(w_o)>1/2$, if $\theta_\eta(w_o)=0$ then $\lambda_{w_o}$ is an eigenvalue for the Friedrichs extension $\widetilde{T}_{\theta_\eta}$.   This extension is self-adjoint and $\widetilde T_{\theta_\eta} \geq 1/4$, so  $\lambda_{w_o}\geq 1/4$ and  $\lambda_{w_o}\in \R$.  Thus when $\theta_\eta(w_o)=0$ and $\text{Re}(w_o)>1/2$, $\lambda_{w_o}$ is an eigenvalue of $\widetilde{T}_{\theta_\eta}$ and $\lambda_{w_o}\geq 1/4$ and  $\lambda_{w_o}\in \R$. No such $w_o$ exist and so all zeros must be on the critical line. 
This proves the main part of Theorem \ref{nozeros}.  In the next section we will show that these zeros are simple.

\section{Simple Zeros}

\tab In order to show that the zeros are simple we will first want to examine $\widetilde\theta_\eta(u_{w})$ for $\text{Re}(w)=1/2$.

{\lem\label{=1/2} For $Re(w_o)=1/2$ and $w_o\neq 1/2$, if $\theta_\eta(w_o)=0$ then $\widetilde\theta_\eta(u_{w_o})=0$ unless $h(w_o)=0=h(1-w_o)$. }

\begin{proof}For $\text{Re}(w_o)=1/2$ ($w_o\neq 1/2$) and $\theta_{\eta}(w_o)=0$, 

$$\widetilde\theta_{\eta}(u_{w_o})=\int_{(1/2)} \frac{{(1+\eta \cdot c_s)}(1+\eta \cdot c_{1-s})}{\lambda_s-\lambda_{w_o}}\,ds \phantom{weeeeeeeeeeeeeeeeeee}$$

$$\phantom{weeeeeeeeeeeee}=\int_{(1/2)} \frac{(1+\eta \cdot c_{1-s})(1+\eta \cdot c_s) - (1+\eta \cdot c_{1-w_o})(1+\eta  \cdot c_{w_o})}{\lambda_s-\lambda_{w_o}}\,ds
$$
$$\phantom{weeeeeeeeeeeeeeeeeeeeeeeeeeeee}
+ \int_{(1/2)} \frac{ (1+\eta \cdot c_{1-w_o})(1+\eta  \cdot c_{w_o})}{\lambda_s-\lambda_{w_o}}\,ds$$

The function \vspace{-.2cm}
$$w\mapsto
\int_{(1/2)} \frac{(1+\eta \cdot c_{1-s})(1+\eta \cdot c_s) - (1+\eta \cdot c_{1-w})(1+\eta  \cdot c_w)}{\lambda_s-\lambda_{w}}\,ds$$
$$\phantom{weeeeeeeeeeeeeeeeeeeeeeeee}+ \int_{(1/2)} \frac{ (1+\eta \cdot c_{1-w})(1+\eta  \cdot c_w)}{\lambda_s-\lambda_{w}}\,ds
$$ 
 is meromorphic in $w$.  As in the proof of Lemma \ref{BC1}, we can evaluate this for $\text{Re}(w)>1/2$ for residues to get that $\displaystyle\widetilde\theta_{\eta}(u_w)=-4\pi i \cdot \frac{1+\eta  \cdot  c_w}{2w-1}.$ Thus at $w=w_o$, it is $$\displaystyle\widetilde\theta_{\eta}(u_{w_o})=-4\pi i \cdot\frac{1+\eta  \cdot c_{w_o}}{2{w_o}-1}=-4\pi i \cdot\frac{\theta_{\eta}(w_o)}{2w_o-1}=0$$
\end{proof}

\tab Finally we have the following. 

{\lem The zeros $w_o$ of $1+\eta\cdot c_w$ on $1/2 +i\mathbb{R}$ are simple except for possibly $w_o=1/2$.}
\begin{proof} Suppose $w_o$ is a zero of $\theta_{\eta}(w)=1+\eta\cdot c_w$.  To show that $w_o$ is simple, we want to show $\theta_{\eta}(w)'(w_o)=\eta\cdot c_w'(w_o)\neq 0. $
Notice $\displaystyle\widetilde\theta_{\eta}(u_w)= \frac{1+\eta\cdot c_w}{1-2w}$ from the proofs of Lemmas \ref{BC1} and \ref{=1/2} and so \vspace{-.2cm}
$$\widetilde\theta_{\eta}(u_w)'(w_o)=-4\pi i \cdot \frac{(2w_o-1)(\eta\cdot c_w)'(w_o)+2(1 +\eta\cdot c_{w_o})}{(2w_o-1)^2}=-4\pi i \cdot\frac{\eta\cdot c_w'(w_o)}{2w_o-1}$$\vspace{-.4cm} 

since $1+\eta\cdot c_{w_o}=0$. Thus in order to show that  $\theta_{\eta}(w)'(w_o)=\eta\cdot c_w'(w_o)\neq 0$ is suffices to show that $\theta_{\eta}(u_w)'(w_o)\neq 0$.

\tab To see that this is in fact non-vanishing, recall that on $\text{Re}(s)=1/2$, 
$$\widetilde\theta_{\eta}(u_w)=\int_{(1/2)} \frac{{(1+\eta\cdot c_s)}(1+\eta\cdot c_{1-s}) - {(1+\eta\cdot c_w)}(1+\eta\cdot c_{1-w})}{\lambda_s-\lambda_{w}}\,ds \phantom{weeeee}(\text{Re}(s)=1/2)$$
$$+ \int_{(1/2)} \frac{{(1+\eta\cdot c_w)}(1+\eta\cdot c_{1-w})}{\lambda_s-\lambda_{w}}\,ds.$$

Taking the derivative in $w$, the last two terms cancel one another, giving

$$\widetilde\theta_{\eta}(u_w)'(w_o)=-\int_{(1/2)} \frac{{(1+\eta\cdot c_s)}(1+\eta\cdot c_{1-s})}{(\lambda_s-\lambda_{w_o})^2}\,ds\cdot (2w_o-1)$$

 Now $\displaystyle\widetilde\theta_{\eta}(u_w)'(w_o)=-(2w_o-1)|u_{w_o}|_{L^2(1/2+i\R)}^2$ which is not zero except at $w_o=1/2$.
\end{proof}

\section{Appendix: Friedrichs Extensions}\label{FE}

\tab  For the convenience of the reader, we will recall some facts about Friedrichs extensions.  The following is given in a general setting which can easily be translated to our case. 

\tab A symmetric, densely-defined operator $S$ on a Hilbert space $V$ is {\it semi-bounded} when $\langle Sv,v\rangle\geq c\cdot\langle v, v\rangle$ or $\langle Sv,v\rangle\leq c\cdot\langle v, v\rangle$ for some real constant $c$. We can construct the Friedrichs extension of a densely-defined, symmetric semi-bounded operator $S$ as follows:

\tab Without loss of generality, consider a densely-defined, symmetric operator $S$  with dense domain $D_S$ and $\langle S v, v\rangle \geq \langle v,v\rangle$  for all $v\in D_S$. 

\tab Define an inner product $\langle, \rangle_1$ on $D_S$ by $\langle v, w\rangle_1:=\langle Sv, w\rangle $ for  $v, w\in D_S$ and let $V^1$ be the completion of $D_S$ with respect to the metric induced by $\langle,\rangle_1$.  Since $\langle v,v\rangle_1\geq\langle v, v\rangle$, the inclusion map $D_S\hookrightarrow V$ extends to a continuous map $V^1\hookrightarrow V$. Furthermore,  $V^1$ is also dense in $V$ since $D_S$ is dense in $V$.
 
 \tab For $w\in V$, the functional $v\mapsto\langle v, w\rangle $ is a continuous linear functional on $V^1$ with norm
 $$\sup_{|v|_1\leq 1}|\langle v, w\rangle|\leq\sup_{|v|_1\leq 1}|v|\cdot|w|\leq\sup_{|v|_1\leq 1}|v|_1\cdot|w|=|w|.$$
 By the Riesz-Fr\'{e}chet Theorem on $V^1$, there is a $w'\in V^1$ so that $\langle v,w'\rangle_1=\langle v, w\rangle$ for all $v\in V^1$ and $w\in V$ with norm bounded by the norm of $v\mapsto\langle v, w\rangle$; explicitly, $|w'|_1\leq|w|$.  The map $A:V\to V^1$ defined by $w\mapsto w'$ is linear.  The densely-defined inverse of $A$ will be a self-adjoint extension $\widetilde S$, the {\it Friedrichs extension} of $S$. The Friedrichs extension $\widetilde S$ is self-adjoint and an extension of $S$.  This is due to Friedrichs \cite{Friedrichs1934} and is also on p.103 of vonNeumann's 1929 paper \cite{vonNeumann1929}. \\
 
 {\thm $\widetilde S$ is a self-adjoint extension of $S$.}\\

\tab This construction proves following theorem of Friedrichs \cite{Friedrichs1934}.\\

{\thm A positive, densely-defined, symmetric operator $S$ with domain $D_S$ has a positive self-adjoint extension with the same lower bound.  }\\

\tab This extension has useful properties of particular interest here.  An alternative characterization of the extension makes this clearer.

\tab Assume that $V$ has a $\mathbb{C}$-linear complex conjugation $v\to v^c$ with the properties: $(v^c)^c=v$ and  $\langle v^c,w^c\rangle=\overline{\langle v, w\rangle}$.  Further, let $S$ commute with conjugation so that $(Sv)^c=S(v^c)$. Let $V^{-1}$ be the dual of $V^1$ so that $V^1\subset V\subset V^{-1}$ via the embedding $v\mapsto (*\mapsto\langle v, *^c\rangle)\in V^{-1}$ for $v\in V$.

\tab Given this small adaptation, there is an alternate characterization of the Friedrichs extension.  To give it, define a continuous, complex-linear map $S^\#:V^1\to V^{-1}$ by 
 $$(S^\#v)(w)=\langle v, w^c\rangle_1$$ for $v, w\in V^1$.\\

  {\thm Let $X=\{v\in V^1~|~S^\#v\in V\}$.  Then the Friedrichs extension of $S$ is $\widetilde S=S^\#|_X$ with domain $D_{\widetilde S}=X$.}
  
  \begin{proof} Let $B=S^\#|_X$.  Let $A:V\to V^1$ be the inverse of $\widetilde S$ defined by $\langle Av, w\rangle_1=\langle v, w\rangle$ for all $w\in V^1$ and $v\in V$ from the Riesz-Fr\'{e}chet Theorem.  Then 
  $$\langle BAv, w\rangle =\langle Av, w\rangle_1=\langle v, w\rangle$$ for $v\in V$ and $w\in V^1$.  Also,   
    $$\langle ABv, w\rangle_1 =\langle Bv, w\rangle=\langle v, w\rangle_1$$ for $v\in X$ and $w\in V^1$.  This $B=A^{-1}=\widetilde S$. 
    \end{proof}
    \vspace{.3cm}

\subsection{Extensions of Restrictions}

Using the latter characterization of the Friedrichs extension we can see how the construction of the extension behaves for restricted operators. Again, the following is given in a general setting which can easily be translated to our case. We will assume that $S$ and the related terms are as defined in Section \ref{FE}.  From above we have the following inclusions $i$ and $i^*$
$$\begin{tikzcd}
V^1\arrow[r, "i"] & V  \arrow[r, "\cong"] &V^* \arrow[r, "i^*"] & V^{-1}
\end{tikzcd}$$
where $V^*$ is the complex linear dual of $V$.

\tab Let $\theta \in V^{-1}$ and assume that $\theta\notin i^*(V^*)$. Note that $\text{ker}\,\theta$ is a closed subspace of $V^1$.  The following lemma follows from the general fact that for a continuous inclusion of Hilbert spaces 
 $i: V^1\to V$ for $D_S\subset V^1$ is dense in $V$ and for a finite-dimensional (in our case, {\it one}-dimensional subspace spanned by $\theta$) subspace $\theta\subset (V^1)^*=V^{-1}$  such that $i^*(V^*)\cap \theta= \{0\}$, we have that $D_S\cap \text{ker}\,\theta\subset V^1$ is dense in $V$.\\

{\lem\label{ker} If $\theta \notin i^*(V^*)$ then $D_S\cap \text{ker}\,\theta  $ is dense in $V$.}

\begin{proof}   Since $\theta\notin i^*(V^*)$ and $\theta  $ cannot be in the $V$-topology on dense $D_S$.  This gives us that there is a $\delta>0$ so that for each $\epsilon>0$ an element $x_{\epsilon}\in D_S$ with $|x_{\epsilon}|_V<\epsilon$ and $|\theta  (x_{\epsilon})|\geq \delta$.  Given $y\in V$ density of $D_S$ in $V$ gives a sequence $z_n$ in $D_S$ approaching $y$ in the $V$-topology.  If $\theta(z_n)=0$ for infinitely many $n$ then we are done.  Otherwise, the define the sequence
 $$z_n'=z_n-\frac{\theta(z_n)}{\theta(x_{\epsilon_n})}x_{\epsilon_n} \phantom{eee}\text{ with } \phantom{eee}\epsilon_n = \frac{\delta}{|\theta(z_n)|}2^{-n}$$ 
 is in $\text{ker}\,\theta  $.  Then 
 $$\theta  (z_n')=\theta(z_n) - \frac{\theta  (z_n)}{\theta  (x_{\epsilon_n})}\theta  (x_{\epsilon_n})=0$$
and $z_n'\to y$ in the $V$-topology since
$$ \left|\frac{\theta  (z_n)}{\theta  (x_{\epsilon_n})}x_{\epsilon_n}\right|_V
=\left|\frac{\theta  (z_n)}{\theta  (x_{\epsilon_n})}\right|\cdot \left|x_{\epsilon_n}\right|_V
<|\theta  (z_n)|\frac{1}{\delta}\cdot \frac{\delta}{|\theta  (z_n)|}2^{-n}=2^{-n}\to 0$$

 \end{proof}

\tab Induction can be used to generalize the result for any such finite-dimensional subspace $\theta  $ of $V^{-1}$.

\tab Define $S_{\theta}:=S|_{D_S\cap\text{ker}\,\theta  }$ then $D_{S_\theta}:=D_S\cap\text{ker}\,\theta$.   Since $\theta  \notin V$ as in Lemma \ref{ker}, $D_{S_\theta}$ is still dense in $V$ and since $S_{\theta}$ is a restriction of $S$, the symmetry and $\langle S_{\theta}\,v,v\rangle\geq\langle v, v\rangle$ properties are inherited from $S$. The $V^1$-closure of $D_{S_\theta}$ is $V^1\cap \text{ker}\,\theta$.

\tab Let $W^{-1}$ be the dual of $W^1=\text{ker}\,\theta$ (on $V^1$) so we have $W^1= \text{ker}\,\theta\subset V\subset W^{-1}$.  This yields the following diagram 
$$\begin{tikzcd}
V^1\arrow[r] & V \arrow[r] & V^{-1}\arrow[d, "j^*"]  \\
W^1=\text{ker}\,{\theta  } \arrow[u]  & & W^{-1}\end{tikzcd}$$

Recall $S^\#:V^1\to V^{-1}$ by 
  $(S^{\#}v)(w):=\langle v, w^c\rangle_1$ for $v,w\in V^1$.  
  
  {\thm\label{Fried} The Friedrichs extension $\widetilde S_{\theta}$ of $S_{\theta}$ has domain $D_{\widetilde S_\theta}=\{v\in W^1~|~S^\#v\in V+\mathbb{C}\cdot \theta\}$ and is characterized by 
$$\widetilde S_{\theta  }\,v=w ~~~\Longleftrightarrow~~~S^\#v\in w+\mathbb{C}\cdot \theta  $$ for $v\in D_{\widetilde S_{\theta}}$ and $w\in  V$.}

\begin{proof} Define $S_{\theta  }^\#:W^1\to W^{-1}$ by $$(S_{\theta}^{\#}v)(w):=\langle v, w^c\rangle_1$$ for all $w\in W^1$.  The domain of the Friedrichs extension $\widetilde S_{\theta  }$ is
$$D_{\widetilde S_{\theta}}=\{v\in W^1~|~ S_{\theta}^\#v\in V\}$$ where $V$ is the cope of $V$ in $W^1\to V\to W^{-1}$ and $\widetilde S_{\theta}=S_{\theta  }^\#|_{D_{\widetilde S_{\theta}}}$ by Theorem 9.   With the inclusion $j:W^1\to V^1$  for all $x, y\in W^1$
$$(S_{\theta  }^\#x)(y)=\langle jx, (jy)^c\rangle_1=(S^\#jx)(jy)=\left( (j^*\circ S^\#\circ j) x\right)(y)$$ 
and so $S_{\theta  }^\#=j^*\circ S^\#\circ j$ and $$D_{\widetilde S_{\theta  }}=\{v\in W^1~|~j^*(S^\#(jv))\in V\}.$$ Furthermore $S^\#(jv)\in V+\mathbb{C}\cdot \theta   $ and the inclusion map $j$ is redundant. 
The dual $W^{-1}$ of $W^1$ is 
$$W^{-1}=(V^1\cap\text{ker}\,\theta )^*\cong V^{-1}/\mathbb{C}\cdot\theta.$$\end{proof}

\tab The Friedrichs extension makes the following diagram commute: 

$$\begin{tikzcd}
V^1\arrow[rr, bend left, "S^\#"]\arrow[r] & V & V^{-1}\arrow[d, "j^*"]  \\
W^1 \arrow[u, "j"]   \arrow[rr, bend left, "S_{\theta}^\#"]\arrow[r] & V \arrow[r] & W^{-1}\\
D_{\tilde S_{\theta}}\arrow[u, hook]\arrow[ru, "\widetilde S_{\theta}"]& &
\end{tikzcd}$$

(The apparent missing arrows are excluded because the diagram would not otherwise commute.)

\section*{Acknowledgements}
The author would like to thank Paul Garrett for his guidance and suggestion of the problem and the reviewer for their helpful feedback. The author also acknowledges support from NSF Grant number DMS-2001909.


\end{document}